\documentclass[leqno]{article}

\def\im{\mathop{\rm Im}}
\def\re{\mathop{\rm Re}}

\newtheorem{theorem}{Theorem}
\newtheorem{lemma}[theorem]{Lemma}
\newtheorem{proposition}[theorem]{Proposition}
\newtheorem{definition}[theorem]{Definition}
\newtheorem{corollary}[theorem]{Corollary}

\newcommand{\begintheorem}{\addtocounter{equation}{1}\begin{theorem}}
\newcommand{\beginlemma}{\addtocounter{equation}{1}\begin{lemma}}
\newcommand{\beginproposition}{\addtocounter{equation}{1}\begin{proposition}}
\newcommand{\begindefinition}{\addtocounter{equation}{1}\begin{definition}}
\newcommand{\begincorollary}{\addtocounter{equation}{1}\begin{corollary}}

\begin{document}

\title{Elements of harmonic analysis, 4}

\author{Stephen William Semmes	\\
	Rice University		\\
	Houston, Texas}

\date{}

\maketitle

\tableofcontents

\bigskip

        These informal notes are based on a course given at Rice
University in the spring semester of 2004, and much more information
can be found in the references.

\section{Gauss--Weierstrass kernels}
\label{gauss--weierstrass kernels}
\setcounter{equation}{0}

	If $z = x + i \, y$ is a complex number, with $x$, $y$ real
numbers, then the complex conjugate of $z$ is denoted $\overline{z}$
and defined to be $x - i \, y$.  We call $x$, $y$ the real and
imaginary parts of $z$, and denote them $\re z$, $\im z$.  Notice that
the complex conjugate of a sum or product of two complex numbers is
equal to the corresponding sum or product of complex conjugates.

	If $x$ is a real number, then the absolute value of $x$ is
denoted $|x|$ and defined to be equal to $x$ when $x \ge 0$ and to
$-x$ when $x \le 0$.  For a complex number $z$ the absolute value or
modulus $|z|$ is defined by $|z| = \sqrt{(\re z)^2 + (\im z)^2}$,
which is the same as saying that $|z|$ is a nonnegative real number
and that $|z|^2 = z \, \overline{z}$.  One can check that the modulus
of a product of complex numbers is equal to the modulus of the
product.  Also, the triangle inequality states that the modulus of a
sum of complex numbers is less than or equal to the sum of the moduli
of the individual complex numbers.

	Fix a positive integer $n$, and let ${\bf R}^n$, ${\bf C}^n$
denote the usual real and complex vector spaces of $n$-tuples of real
and complex numbers, with respect to coordinatewise addition and
scalar multiplication.  If $z = (z_1, \ldots, z_n)$, $w = (w_1,
\ldots, w_n)$ are elements of ${\bf C}^n$, we put
\begin{equation}
	z \cdot w = \sum_{j=1}^n z_j \, w_j
\end{equation}
and
\begin{equation}
	|z| = \biggl(\sum_{j=1}^n |z_j|^2\biggr)^{1/2}.
\end{equation}
Thus
\begin{equation}
	|x|^2 = x \cdot x
\end{equation}
when $x \in {\bf R}^n$.  For $z, w \in {\bf C}^n$ we have the
Cauchy--Schwarz inequality
\begin{equation}
	|z \cdot w| \le |z| \, |w|
\end{equation}
and the triangle inequality
\begin{equation}
	|z + w| \le |z| + |w|.
\end{equation}

	For each positive real number define the corresponding
Gauss and Weierstrass kernels by
\begin{equation}
	G_\alpha(x) = \exp (-4\pi^2 \, \alpha \, x \cdot x)
\end{equation}
and
\begin{equation}
	W_\alpha(x) = (4 \pi \, \alpha)^{-n/2} \exp (- x \cdot x / (4 \alpha)).
\end{equation}
We might normally think of these as functions on ${\bf R}^n$, and they
extend to complex-analytic functions on ${\bf C}^n$ in an obvious
manner as well, just using the same formulae.  Of course they are also
basically the same functions, except for some modest adjustments.

	Two important identities involving these functions are
\begin{equation}
	W_\alpha(\xi) 
   = \int_{{\bf R}^n} G_\alpha(x) \, \exp(- 2 \pi i \, x \cdot \xi) \, dx
\end{equation}
and
\begin{equation}
	G_\alpha(\xi) 
   = \int_{{\bf R}^n} W_\alpha(x) \, \exp(- 2 \pi i \, x \cdot \xi) \, dx.
\end{equation}
These identities are essentially the same as each other, and they are
also essentially the same for any choice of $\alpha > 0$, simply by
rescaling.

	When $\xi = 0$ these identities reduce to the well-known fact that
\begin{equation}
	\int_{-\infty}^\infty \exp (- \pi \, x^2) \, dx = 1.
\end{equation}
Namely, the integral here is obviously a positive real number, and it
suffices to show that its square is equal to $1$.  To do this one can
identify its square with a $2$-dimensional integral of the same type.
The $2$-dimensional integral can be converted into a $1$-dimensional
integral again using polar coordinates, and the $1$-dimensional
integral can be integrated directly.

	The identities above make sense for all $\xi \in {\bf C}^n$.
If $\xi$ is purely imaginary, which is to say that $\xi$ is equal to
$i$ times an element of ${\bf R}^n$, then one can reduce the integrals
to the $\xi = 0$ case by making a change of variables.  To deal with the
general case one can argue that the two sides of the equations
define complex-analytic functions of $\xi \in {\bf C}^n$ which agree
when $\xi$ is purely imaginary, and hence agree for all $\xi$.

\section{Continuous functions on ${\bf R}^n$}
\label{continuous functions on R^n}
\setcounter{equation}{0}

	Let $f(x)$ be a continuous complex-valued function on ${\bf
R}^n$.  We say that $f$ is integrable on ${\bf R}^n$ if $\int_{{\bf
R}^n} |f(x)| \, dx < \infty$, where the integral can be defined
through standard methods as an improper integral.  In this event
$\int_{{\bf R}^n} f(x) \, dx$ can also be defined as an improper
integral, and in a rather stable way, so that various standard
approximations can be employed and yield the same result.

	Now suppose that $f(x)$ is a complex-valued continuous function
on ${\bf R}^n$ which is, say, bounded.  One way to approach the integral
of $f(x)$ is to use Gauss means, which is to say the integrals
\begin{equation}
	\int_{{\bf R}^n} f(x) \, G_\alpha(x) \, dx
\end{equation}
for $\alpha > 0$.  These integrals make sense because the Gaussian
functions are integrable, and if the limit of these integrals as
$\alpha \to 0$ exists, then we say that $f(x)$ is Gauss summable to
the limit.  If $f(x)$ is a nonnegative real-valued function, then
this is equivalent to saying that $f(x)$ is integrable.

	If $f(x)$ is a complex-valued continuous function on ${\bf
R}^n$ which is bounded or integrable, consider the function
\begin{equation}
	(W_\alpha * f)(x) = \int_{{\bf R}^n} f(y) \, W_\alpha(x - y) \, dy.
\end{equation}
This is basically an average of the values of $f$ around $x$, since
$W_\alpha$ is nonnegative and has integral equal to $1$, and one can
check that $W_\alpha * f$ is also continuous.  If $f(x)$ is
integrable, then
\begin{equation}
	\int_{{\bf R}^n} |(W_\alpha * f)(x)| \, dx
		\le \int_{{\bf R}^n} |f(y)| \, dy,
\end{equation}
and if $f(x)$ is bounded, then
\begin{equation}
	\sup \{|(W_\alpha * f)(x)| : x \in {\bf R}^n\}
		\le \sup \{|f(y)| : y \in {\bf R}^n\}.
\end{equation}

	Assuming that $f(x)$ is bounded or integrable, one can show that
\begin{equation}
	\lim_{\alpha \to 0} (W_\alpha * f)(x) = x
\end{equation}
for all $x \in {\bf R}^n$.  Moreover one has uniform convergence for
$x$ in any fixed compact subset of ${\bf R}^n$, basically because $f$
is uniformly continuous on compact subsets of ${\bf R}^n$.  If $f(x)$
is uniformly continuous on all of ${\bf R}^n$, then we have uniform
convergence on all of ${\bf R}^n$.

\section{Applications}
\label{applications}
\setcounter{equation}{0}

	Let $f(x)$ be a continuous complex-valued function on ${\bf R}^n$
which is integrable.  Define the Fourier transform of $f(x)$ by
\begin{equation}
	\widehat{f}(\xi) 
	   = \int_{{\bf R}^n} f(x) \, \exp(- 2 \pi i \, x \cdot \xi) \, dx.
\end{equation}
Thus $G_\alpha$ and $W_\alpha$ are Fourier transforms of each other.
In general one can show that $\widehat{f}(\xi)$ is uniformly
continuous when $f$ is integrable, and also that $\widehat{f}(\xi)$ is
bounded, with
\begin{equation}
	\sup \{|\widehat{f}(\xi)| : \xi \in {\bf R}^n \}
		\le \int_{{\bf R}^n} |f(x)| \, dx.
\end{equation}

	Now suppose that $\phi$ is an integrable complex-valued
continuous function on ${\bf R}^n$, and define its inverse Fourier
transform by
\begin{equation}
	{\check \phi}(x) 
	  = \int_{{\bf R}^n} \phi(\xi) \, \exp(2 \pi i \, x \cdot \xi) \, d\xi.
\end{equation}
Just as for the Fourier transform, $G_\alpha$ and $W_\alpha$ are
inverse Fourier transforms of each other for all $\alpha > 0$.  Also,
one can show that $\check \phi$ is a bounded uniformly continuous
function on ${\bf R}^n$ when $\phi$ is integrable, with $|{\check
\phi}(x)|$ less than or equal to $\int_{{\bf R}^n} |\phi(\xi)| \,
d\xi$.

	As the names suggest, one would like to show that the inverse
Fourier transform of the Fourier transform of $f$ is equal to $f$ when
$f(x)$ is a continuous complex-valued function on ${\bf R}^n$ which is
integrable.  However, we do not know in general that the Fourier
transform of $f$ is integrable in this case.  We do know that it is
bounded, and so we can try to look at the inverse Fourier transform
using Gauss summability.

	Thus let $f(x)$ be an integrable complex-valued continuous
function on ${\bf R}^n$, and consider the integral
\begin{equation}
\label{gauss sum of inverse fourier transform of f}
	\int_{{\bf R}^n} \widehat{f}(\xi) \, \exp (2 \pi i \, x \cdot \xi) \,
				G_\alpha(\xi) \, d\xi
\end{equation}
for $\alpha > 0$.  Let us write this as $\int_{{\bf R}^n}
\widehat{f}(\xi) \, \psi(\xi) \, d\xi$ with $\psi(\xi) = G_\alpha(\xi)
\, \exp (2 \pi i x \cdot \xi)$.  It is not too difficult to verify
that this is equal to $\int_{{\bf R}^n} f(y) \, \widehat{\psi}(y) \,
dy$.  This is the ``Multiplication Formula'' for the Fourier
transform.

	For any continuous complex-valued function $h(u)$ on ${\bf
R}^n$ which is integrable, we have that $h(u) \exp(2 \pi i \, a \cdot
u)$ is integrable and the Fourier transform of it at $\eta$ is equal
to $\widehat{h}(\eta - a)$.  The Fourier transform of $G_\alpha$ is
equal to $W_\alpha$, and thus $\widehat{\psi}(y) = W_\alpha(y - x)$ in
the context of the preceding paragraph.  Therefore (\ref{gauss sum of
inverse fourier transform of f}) is equal to $(W_\alpha * f)(x)$, and
it follows that the inverse Fourier transform of $\widehat{f}$ is
Gauss summable to $f(x)$.

	By a bounded measure on ${\bf R}^n$ let us mean a linear
functional $\lambda$ from the vector space of bounded continuous
complex-valued functions on ${\bf R}^n$ into the complex numbers
which is bounded, in the sense that
\begin{equation}
	|\lambda(h)| \le L \, \sup \{|h(x)| : x \in {\bf R}^n\}
\end{equation}
for some positive real number $L$ and all bounded continuous functions
$h$ on ${\bf R}^n$, and which satisfies the continuity property that
if $\{h_j\}_{j=1}^\infty$ is a sequence of bounded continuous
functions on ${\bf R}^n$ which are uniformly bounded and converge to
another bounded continuous function $h$ on ${\bf R}^n$ uniformly on
compact subsets of ${\bf R}^n$, then
\begin{equation}
	\lim_{j \to \infty} \lambda(h_j) = \lambda(h).
\end{equation}
As a basic class of examples, if $f(x)$ is a continuous complex-valued
function on ${\bf R}^n$ which is integrable, then
\begin{equation}
	\lambda_f(h) = \int_{{\bf R}^n} f(x) \, h(x) \, dx
\end{equation}
defines a bounded measure on ${\bf R}^n$, with $L$ equal to
$\int_{{\bf R}^n} |f(x)| \, dx$.  As another class of examples, if $a
\in {\bf R}^n$, then we get a bounded measure on ${\bf R}^n$, the
Dirac mass at $a$, defined by evaluating a given function at $a$.

	If $\lambda$ is a bounded measure on ${\bf R}^n$, then
we define its Fourier transform $\widehat{\lambda}$ by
\begin{equation}
	\widehat{\lambda}(\xi) = \lambda(e_\xi),
		\quad e_\xi(x) = \exp (- 2 \pi i \, x \cdot \xi).
\end{equation}
One can check that this defines a bounded continuous function on ${\bf
R}^n$, and of course this reduces to the previous definition when
$\lambda = \lambda_f$ for an integrable function $f$.  One can also
define $(W_\alpha * \lambda)(y)$ by applying $\lambda$ to $W_\alpha(y
- x)$, viewed as a function of $y$ with $x$ fixed.

	The inverse Fourier transform of $\widehat{\lambda}$ can be
analyzed again using the Gauss means
\begin{equation}
	\int_{{\bf R}^n} \widehat{\lambda}(\xi)
		\, \exp(2 \pi i \, x \cdot \xi) \, G_\alpha(\xi) \, d\xi.
\end{equation}
As in the earlier calculations this reduces to $(W_\alpha * \lambda)(x)$.
This converges to $\lambda$ as $\alpha \to 0$ in the weak sense that
\begin{equation}
	\int_{{\bf R}^n} (W_\alpha * \lambda)(x) \, h(x) \, dx
\end{equation}
converges to $\lambda(h)$ as $\alpha \to 0$ for all bounded continuous
functions $h$ on ${\bf R}^n$, basically because the preceding integral
is the same as $\lambda$ applied to $W_\alpha * h$, and $W_\alpha * h$
is uniformly bounded and converges to $h$ uniformly on compact subsets
of ${\bf R}^n$ as $\alpha \to 0$.

\end{document}